\documentstyle[12pt,fleqn]{article}

\begin{document}

\title{Complex Numbers in 5 Dimensions}

\author{Silviu Olariu
\thanks{e-mail: olariu@ifin.nipne.ro}\\
Institute of Physics and Nuclear Engineering, Tandem Laboratory\\
76900 Magurele, P.O. Box MG-6, Bucharest, Romania}

\date{4 August 2000}

\maketitle

\abstract

A system of commutative complex numbers in 5 dimensions of the form
$u=x_0+h_1x_1+h_2x_2+h_3x_3+h_4x_4$ is described in this paper, the variables
$x_0, x_1, x_2, x_3, x_4$ being real numbers.  The operations of addition and
multiplication of the 5-complex numbers introduced in this work have a
geometric interpretation based on the the modulus $d$, the amplitude $\rho$,
the polar angle $\theta_+$, the planar angle $\psi_1$, and the azimuthal angles
$\phi_1,\phi_2$.  The exponential function of a 5-complex number can be
expanded in terms of polar 5-dimensional cosexponential functions $g_{5k}(y),
k=0,1,2,3,4$, and the expressions of these functions are obtained from the
properties of the exponential function of a 5-complex variable.  Exponential
and trigonometric forms are obtained for the 5-complex numbers, which depend on
the modulus, the amplitude and the angular variables.  The 5-complex functions
defined by series of powers are analytic, and the partial derivatives of the
components of the 5-complex functions are closely related.  The integrals of
5-complex functions are independent of path in regions where the functions are
regular. The fact that the exponential form of the 5-complex numbers depends on
the cyclic variables $\phi_1, \phi_2$ leads to the concept of pole and residue
for integrals on closed paths. The polynomials of 5-complex variables can be
written as products of linear or quadratic factors.

\endabstract

\section{Introduction}

A regular, two-dimensional complex number $x+iy$ 
can be represented geometrically by the modulus $\rho=(x^2+y^2)^{1/2}$ and 
by the polar angle $\theta=\arctan(y/x)$. The modulus $\rho$ is multiplicative
and the polar angle $\theta$ is additive upon the multiplication of ordinary 
complex numbers.

The quaternions of Hamilton are a system of hypercomplex numbers
defined in four dimensions, the
multiplication being a noncommutative operation, \cite{1} 
and many other hypercomplex systems are
possible, \cite{2a}-\cite{2b} but these hypercomplex systems 
do not have all the required properties of regular, 
two-dimensional complex numbers which rendered possible the development of the 
theory of functions of a complex variable.

A system of complex numbers in 5 dimensions is described in this work,
for which the multiplication is associative and commutative, and which is 
rich enough in properties so that an exponential form exists and the concepts
of analytic 5-complex 
function,  contour integration and residue can be defined.
The 5-complex numbers introduced in this work have 
the form $u=x_0+h_1x_1+h_2x_2+h_3x_3+h_4x_4$, the variables 
$x_0,x_1,x_2,x_3,x_4$ being real numbers. 
If the 5-complex number $u$ is
represented by the point $A$ of coordinates $x_0,x_1,x_2,x_3,x_4$, 
the position of the point $A$ can be described
by the modulus $d=(x_0^2+x_1^2+x_2^2+x_3^2+x_4^2)^{1/2}$, 
by 2 azimuthal angles $\phi_1, \phi_2$, by 1 planar angle $\psi_1$,
and by 1 polar angle $\theta_+$. 

The exponential function of a 5-complex number can be expanded in terms of
the polar 5-dimensional cosexponential functions
$g_{5k}(y)$, $k=0,1,2,3,4$. The expressions of these functions are obtained
from the properties of the exponential function of a 5-complex variable.
Addition theorems and other relations
are obtained for the polar 5-dimensional cosexponential functions.
Exponential and trigonometric forms are given for the 5-complex numbers.
Expressions are obtained for the elementary functions of 5-complex variable.
The functions $f(u)$ of 5-complex variable which are defined by power series
have derivatives independent of the direction of approach to the point under
consideration. 
If the 5-complex function $f(u)$ 
of the 5-complex variable $u$ is written in terms of 
the real functions $P_k(x_0,x_1,x_2,x_3,x_4), k=0,1,2,3,4$, then
relations of equality  
exist between partial derivatives of the functions $P_k$. 
The integral $\int_A^B f(u) du$ of a 5-complex
function between two points $A,B$ is independent of the path connecting $A,B$,
in regions where $f$ is regular.
The fact that the exponential form
of the 5-complex numbers depends on the cyclic variables $\phi_1,\phi_2$
leads to the 
concept of pole and residue for integrals on closed paths,
and if $f(u)$ is an analytic 5-complex function, then $\oint_\Gamma
f(u)du/(u-u_0)$ is expressed in this work in terms of the 5-complex residue
$f(u_0)$. The polynomials of
5-complex variables can be written as products of linear or quadratic
factors. 

This paper belongs to a series of studies on commutative complex numbers in $n$
dimensions. \cite{2c}
The 5-complex numbers described in this work are a particular case for 
$n=5$ of the polar complex numbers in $n$ dimensions.\cite{2c},\cite{2d}

\section{Operations with polar complex numbers in 5 dimensions}

A polar hypercomplex number $u$ in 5 dimensions 
is represented as  
\begin{equation}
u=x_0+h_1x_1+h_2x_2+h_3x_3+h_4x_4. 
\label{1a}
\end{equation}
The multiplication rules for the bases 
$h_1, h_2, h_3, h_4$ are 
\begin{eqnarray}
\lefteqn{h_1^2=h_2,\;h_2^2=h_4,\; h_3^2=h_1,\; h_4^2=h_3,\nonumber}\\
&& h_1h_2=h_3, \;h_1h_3=h_4,\;
h_1h_4=1,\; h_2h_3=1,\; h_2h_4=h_1, h_3h_4=h_2.
\label{1}
\end{eqnarray}
The significance of the composition laws in Eq.
(\ref{1}) can be understood by representing the bases $h_j, h_k$ by points on a
circle at the angles $\alpha_j=2\pi j/5,\alpha_k=2\pi k/5$, as shown in Fig. 1,
and the product $h_j h_k$ by the point of the circle at the angle 
$2\pi (j+k)/5$. If $2\pi\leq 2\pi (j+k)/5<4\pi$, the point represents the basis
$h_l$ of angle $\alpha_l=2\pi(j+k)/5-2\pi$.

The sum of the 5-complex numbers $u$ and
$u^\prime$ is
\begin{equation}
u+u^\prime=x_0+x^\prime_0+h_1(x_1+x^\prime_1)
+h_2(x_2+x^\prime_2)+h_3(x_3+x^\prime_3)+h_4(x_4+x^\prime_4).
\label{2}
\end{equation}
The product of the numbers $u, u^\prime$ is then
\begin{equation}
\begin{array}{l}
uu^\prime=x_0 x_0^\prime +x_1x_4^\prime+x_2 x_3^\prime+x_3x_2^\prime
+x_4x_1^\prime\\
+h_1(x_0 x_1^\prime+x_1x_0^\prime+x_2x_4^\prime+x_3x_3^\prime
+x_4 x_2^\prime) \\
+h_2(x_0 x_2^\prime+x_1x_1^\prime+x_2x_0^\prime+x_3x_4^\prime
+x_4 x_3^\prime) \\
+h_3(x_0 x_3^\prime+x_1x_2^\prime+x_2x_1^\prime+x_3x_0^\prime
+x_4 x_4^\prime) \\
+h_4(x_0 x_4^\prime+x_1x_3^\prime+x_2x_2^\prime+x_3x_1^\prime
+x_4 x_0^\prime). \\
\end{array}
\label{3}
\end{equation}

The relation between the variables
$v_+, v_1, \tilde v_1, v_2, \tilde v_2$ and $x_0, x_1,
x_2, x_3, x_4$ can be written with the aid of the parameters $p=(\sqrt{5}-1)/4,
q=\sqrt{(5+\sqrt{5})/8}$ as
\begin{equation}
\left(
\begin{array}{c}
v_+\\
v_1\\
\tilde v_1\\
v_2\\
\tilde v_2\\
\end{array}\right)
=\left(
\begin{array}{ccccc}
1&1&1&1&1\\
1&p&2p^2-1&2p^2-1&p\\
0&q&2pq&-2pq&-q\\
1&2p^2-1&p&p&2p^2-1\\
0&2pq&-q&q&-2pq\\
\end{array}
\right)
\left(
\begin{array}{c}
x_0\\
x_1\\
x_2\\
x_3\\
x_4
\end{array}
\right).
\label{9e}
\end{equation}
The other variables are $v_3=v_2, \tilde v_3=-\tilde v_2, v_4=v_1, 
\tilde v_4=-\tilde v_1$. 
The variables $v_+, v_1, \tilde v_1, v_2, \tilde v_2$ will be called canonical
5-complex variables.

\section{Geometric representation of polar complex numbers in 5
dimensions}

The 5-complex number $x_0+h_1x_1+h_2x_2+h_3x_3+h_4x_4$
can be represented by 
the point $A$ of coordinates $(x_0,x_1,x_2, x_3, x_4)$. 
If $O$ is the origin of the 5-dimensional space,  the
distance from the origin $O$ to the point $A$ of coordinates
$(x_0,x_1,x_2, x_3, x_4)$ has the expression
\begin{equation}
d^2=x_0^2+x_1^2+x_2^2+x_3^2+x_4^2.
\label{10}
\end{equation}
The quantity $d$ will be called modulus of the 5-complex number 
$u$. The modulus of a 5-complex number $u$ will be designated by $d=|u|$.
The modulus has the property that
\begin{equation}
|u^\prime u^{\prime\prime}|\leq \sqrt{5}|u^\prime||u^{\prime\prime}| .
\label{79}
\end{equation}

The exponential and trigonometric forms of the 5-complex number $u$ can be
obtained conveniently in a rotated system of axes defined by the transformation
\begin{equation}
\left(
\begin{array}{c}
\xi_+\\
\xi_1\\
\eta_1\\
\xi_2\\
\eta_2\\
\end{array}\right)
=\sqrt{\frac{2}{5}}\left(
\begin{array}{ccccc}
\frac{1}{\sqrt{2}}&\frac{1}{\sqrt{2}}
&\frac{1}{\sqrt{2}}&\frac{1}{\sqrt{2}}&\frac{1}{\sqrt{2}}\\
1&p&2p^2-1&2p^2-1&p\\
0&q&2pq&-2pq&-q\\
1&2p^2-1&p&p&2p^2-1\\
0&2pq&-q&q&-2pq\\
\end{array}
\right)
\left(
\begin{array}{c}
x_0\\
x_1\\
x_2\\
x_3\\
x_4
\end{array}
\right).
\label{12}
\end{equation}
The lines of the matrices in Eq. (\ref{12}) gives the components
of the 5 basis vectors of the new system of axes. These vectors have unit
length and are orthogonal to each other.
The relations between the two sets of variables are
\begin{equation}
v_+= \sqrt{5}\xi_+ ,  
v_k= \sqrt{\frac{5}{2}}\xi_k , \tilde v_k= \sqrt{\frac{5}{2}}\eta_k, k=1,2 .
\label{12b}
\end{equation}

The radius $\rho_k$ and the azimuthal angle $\phi_k$ in the plane of the axes
$v_k,\tilde v_k$  are
\begin{eqnarray}
\rho_k^2=v_k^2+\tilde v_k^2, \:\cos\phi_k=v_k/\rho_k,
\:\sin\phi_k=\tilde v_k/\rho_k, 
\label{19a}
\end{eqnarray}
$0\leq \phi_k<2\pi,\; k=1,2$,
so that there are 2 azimuthal angles.
The planar angle $\psi_1$ is
\begin{equation}
\tan\psi_1=\rho_1/\rho_2, 
\label{19b}
\end{equation}
where $0\leq\psi_1\leq\pi/2$.
There is a polar angle $\theta_+$, 
\begin{equation}
\tan\theta_+=\frac{\sqrt{2}\rho_1}{v_+}, 
\label{19c}
\end{equation}
where $0\leq\theta_+\leq\pi$.
It can be checked that
\begin{equation}
\frac{1}{5}v_+^2
+\frac{2}{5}(\rho_1^2+\rho_2^2)=d^2 .
\label{18}
\end{equation}
The amplitude of a 5-complex number $u$ is
\begin{equation}
\rho=\left(v_+\rho_1^2 \rho_2^2\right)^{1/5}.
\label{50bb}
\end{equation}

If $u=u^\prime u^{\prime\prime}$, the parameters of the hypercomplex numbers
are related by
\begin{equation}
v_+=v_+^\prime v_+^{\prime\prime},  
\label{21a}
\end{equation}
\begin{equation}
\rho_k=\rho_k^\prime\rho_k^{\prime\prime}, 
\label{21b}
\end{equation}
\begin{equation}
\tan\theta_+=\frac{1}{\sqrt{2}}\tan\theta_+^\prime \tan\theta_+^{\prime\prime},
\label{21c}
\end{equation}
\begin{equation}
\tan\psi_1=\tan\psi_1^\prime \tan\psi_1^{\prime\prime},  
\label{21d}
\end{equation}
\begin{equation}
\phi_k=\phi_k^\prime+\phi_k^{\prime\prime},
\label{21e}
\end{equation}
\begin{equation}
v_k=v_k^\prime v_k^{\prime\prime}
-\tilde v_k^\prime \tilde v_k^{\prime\prime},\;
\tilde v_k=v_k^\prime \tilde v_k^{\prime\prime}
+\tilde v_k^\prime v_k^{\prime\prime},
\label{22}
\end{equation}
\begin{equation}
\rho=\rho^\prime\rho^{\prime\prime} ,
\label{24}
\end{equation}
where $k=1,2$.

The 5-complex number
$u=x_0+h_1x_1+h_2x_2+h_3x_3+h_4x_4$ can be  represented by the matrix
\begin{equation}
U=\left(
\begin{array}{ccccc}
x_0     &   x_1     &   x_2   &   x_3  &  x_4\\
x_4     &   x_0     &   x_1   &   x_2  &  x_3\\
x_3     &   x_4     &   x_0   &   x_1  &  x_2\\
x_2     &   x_3     &   x_4   &   x_0  &  x_1\\
x_1     &   x_2     &   x_3   &   x_4  &  x_0\\
\end{array}
\right).
\label{24b}
\end{equation}
The product $u=u^\prime u^{\prime\prime}$ is
represented by the matrix multiplication $U=U^\prime U^{\prime\prime}$.

\section{The polar 5-dimensional cosexponential functions}

The polar cosexponential functions in 5 dimensions are
\begin{equation}
g_{5k}(y)=\sum_{p=0}^\infty y^{k+5p}/(k+5p)!, 
\label{29}
\end{equation}
for $k=0,...,4$.
The polar cosexponential functions $g_{5k}$ do not have a definite
parity. 
It can be checked that
\begin{equation}
\sum_{k=0}^{4}g_{5k}(y)=e^y.
\label{29a}
\end{equation}
The exponential of the quantity $h_k y, k=1,...,4$ can be written as
\begin{equation}
\begin{array}{l}
e^{h_1 y}=g_{50}(y)+h_1g_{51}(y)+h_2g_{52}(y)+h_3g_{53}(y)+h_4g_{54}(y),\\
e^{h_2 y}=g_{50}(y)+h_1g_{53}(y)+h_2g_{51}(y)+h_3g_{54}(y)+h_4g_{52}(y),\\
e^{h_3 y}=g_{50}(y)+h_1g_{52}(y)+h_2g_{54}(y)+h_3g_{51}(y)+h_4g_{53}(y),\\
e^{h_4 y}=g_{50}(y)+h_1g_{54}(y)+h_2g_{53}(y)+h_3g_{52}(y)+h_4g_{51}(y).
\end{array}
\label{28b}
\end{equation}

The polar cosexponential functions in 5 dimensions can be obtained by
calculating $e^{(h_1+h_4)y}$ and $e^{(h_1-h_4)y}$ and then by
nultiplying the resulting expression.
The series expansions for $e^{(h_1+h_4)y}$ and $e^{(h_1-h_4)y}$ are
\begin{equation}
e^{(h_1+h_4)y}=\sum_{m=0}^\infty \frac{1}{m!}(h_1+h_4)^m y^m ,
\label{g5-1}
\end{equation}
\begin{equation}
e^{(h_1-h_4)y}=\sum_{m=0}^\infty \frac{1}{m!}(h_1-h_4)^m y^m .
\label{g5-2}
\end{equation}
The powers of $h_1+h_4$ have the form
\begin{equation}
(h_1+h_4)^m=A_m (h_1+h_4)+B_m (h_2+h_3)+C_m.
\label{g5-3}
\end{equation}
The recurrence relations for $A_m, B_m, C_m$ are
\begin{equation}
A_{m+1}=B_m+C_m, B_{m+1}=A_m+B_m, C_{m+1}=2A_m,
\label{g5-4}
\end{equation}
and $A_1=1, B_1=0, C_1=0, A_2=0, B_2=1, C_2=2, A_3=3, B_3=1, C_3=0$.
The expressions of the coefficients are
\begin{equation}
A_m=\frac{2^m}{5}+\frac{2-3a}{5}a^{m-3}+(-1)^{m-3}\frac{5+3a}{5}(1+a)^{m-3},
m\geq 3 ,
\label{g5-5}
\end{equation}
\begin{equation}
B_m=\frac{2^m}{5}+\frac{a-1}{5}a^{m-3}-(-1)^{m-3}\frac{a+2}{5}(1+a)^{m-3},
m\geq 3,
\label{g5-6}
\end{equation}
\begin{equation}
C_m=\frac{2^m}{5}+\frac{4-6a}{5}a^{m-4}+(-1)^{m-4}\frac{10+6a}{5}(1+a)^{m-4} ,
m\geq 4,
\label{g5-7}
\end{equation}
where $a$ is a solution of the equation 
$a^2+a-1=0$.
Substituting the expressions of $A_m, B_m, C_m$ from Eqs.
(\ref{g5-5})-(\ref{g5-7}) in Eq. (\ref{g5-1}) and grouping the terms yields
\begin{eqnarray}
e^{(h_1+h_4)y}
\lefteqn{=\frac{1}{5}e^{2y}+\frac{2}{5}e^{ay}+\frac{2}{5}e^{-(1+a)y}
+(h_1+h_4)
\left[\frac{1}{5}e^{2y}+\frac{a}{5}e^{ay}
-\frac{a+1}{5}e^{-(1+a)y}\right]\nonumber}\\
&&+(h_2+h_3)
\left[\frac{1}{5}e^{2y}-\frac{a+1}{5}e^{ay}+\frac{a}{5}e^{-(1+a)y}\right].
\label{g5-8}
\end{eqnarray}

The odd powers of $h_1-h_4$ have the form
\begin{equation}
(h_1-h_4)^{2m+1}=D_m (h_1-h_4)+E_m (h_2-h_3).
\label{g5-9}
\end{equation}
The recurrence relations for $D_m, E_m$ are
\begin{equation}
D_{m+1}=-3D_m-E_m, E_{m+1}=-D_m-2E_m,
\label{g5-10}
\end{equation}
and $D_1=-3, E_1=-1, D_2=10, E_2=5$.
The expressions of the coefficients are
\begin{equation}
D_m=(b+1)b^{m-1}+(-1)^{m-2}(b+4)(5+b)^{m-1}, m\geq 1,
\label{g5-11}
\end{equation}
\begin{equation}
E_m=-\frac{b+1}{b+2}b^{m-1}+\frac{(-1)^{m-2}}{b+2}(5+b)^{m-1}, 
m\geq 1,
\label{g5-12}
\end{equation}
where $b$ is a solution of the equation 
$b^2+5b+5=0$.
The even powers of $h_1-h_4$ have the form
\begin{equation}
(h_1-h_4)^{2m}=F_m (h_1+h_4)+G_m (h_2+h_3)+H_m.
\label{g5-13}
\end{equation}
The recurrence relations for $F_m, G_m, H_m$ are
\begin{equation}
F_{m+1}=-F_m+G_m, G_{m+1}=F_m-2G_m+H_m, H_{m+1}=2(G_m-H_m),
\label{g5-14}
\end{equation}
and $F_1=0, G_1=1, H_1=-2, F_2=1, G_2=-4, H_2=6$.
The expressions of the coefficients are
\begin{equation}
F_m=-\frac{1}{5(b+2)}b^m+(-1)^{m-1}\frac{b+1}{5(b+2)}(5+b)^m, m\geq 1,
\label{g5-15}
\end{equation}
\begin{equation}
G_m=\frac{4b+5}{5(b+2)}b^{m-1}+(-1)^{m-1}\frac{1}{5(b+2)}(5+b)^m,
m\geq 1,
\label{g5-16}
\end{equation}
\begin{equation}
H_m=-\frac{6b+10}{5(b+2)}b^{m-1}+(-1)^m\frac{2}{5}(5+b)^m,
m\geq 1,
\label{g5-17}
\end{equation}
where $b$ is a solution of the equation 
$b^2+5b+5=0$.

Substituting the expressions of $D_m, E_m, F_m, G_m, H_m$ from Eqs.
(\ref{g5-11})-(\ref{g5-12}) 
and (\ref{g5-15})-(\ref{g5-17})in Eq. (\ref{g5-2}) and grouping the terms
yields 
\begin{eqnarray}
\lefteqn{e^{(h_1-h_4)y}
=\frac{1}{5}+\frac{2}{5}\cos(\sqrt{-b}y)
+\frac{2}{5}\cos(\sqrt{5+b}y)\nonumber}\\
&&+(h_1+h_4)
\left[\frac{1}{5}-\frac{b+3}{5}\cos(\sqrt{-b}y)
+\frac{b+2}{5}\cos(\sqrt{5+b}y)\right]\nonumber\\
&&+(h_2+h_3)
\left[\frac{1}{5}+\frac{b+2}{5}\cos(\sqrt{-b}y)
-\frac{b+3}{5}\cos(\sqrt{5+b}y)\right]\nonumber\\
&&+(h_1-h_4)
\left[\frac{\sqrt{-b}}{5}\sin(\sqrt{-b}y)
+\frac{1}{\sqrt{-5b}}\sin(\sqrt{5+b}y)\right]\nonumber\\
&&+(h_2-h_3)
\left[-\frac{2b+5}{5\sqrt{-b}}\sin(\sqrt{-b}y)
+\frac{b+2}{\sqrt{-5b}}\sin(\sqrt{5+b}y)\right].
\label{g5-18}
\end{eqnarray}

On the other hand, $e^{2h_1 y}$ can be written with the aid of the
5-dimensional polar cosexponential functions as
\begin{equation}
e^{2h_1 y}=g_{50}(2y)+h_1 g_{51}(2y)+h_2 g_{52}(2y)+h_3
g_{53}(2y)+h_4 g_{54}(2y).
\label{g5-19}
\end{equation}
The multiplication of the expressions of $e^{(h_1+h_4)y}$ and
$e^{(h_1-h_4)y}$ in Eqs. (\ref{g5-8}) and (\ref{g5-18}) and the
separation of the real components yields the expressions of the 5-dimensional
cosexponential functions, for $a=(\sqrt{5}-1)/2, b=-(5+\sqrt{5})/2$, as
\begin{eqnarray}
g_{50}(2y)=\frac{1}{5}e^{2y}
+\frac{2}{5}e^{ay}\cos(\sqrt{-b}y)
+\frac{2}{5}e^{-(1+a)y}\cos(\sqrt{5+b}y),
\label{g5-20}
\end{eqnarray}
\begin{eqnarray}
\lefteqn{g_{51}(2y)=\frac{1}{5}e^{2y}
+\frac{1}{5}e^{ay}\left[\frac{-1+\sqrt{5}}{2}
\cos(\sqrt{-b}y)+\frac{5+\sqrt{5}}{2\sqrt{-b}}
\sin(\sqrt{-b}y)\right]\nonumber}\\
&&+\frac{1}{5}e^{-(1+a)y}\left[-\frac{1+\sqrt{5}}{2}
\cos(\sqrt{5+b}y)+\sqrt{\frac{5}{-b}}\sin(\sqrt{5+b}y)\right],
\label{g5-21}
\end{eqnarray}
\begin{eqnarray}
\lefteqn{g_{52}(2y)=\frac{1}{5}e^{2y}
+\frac{1}{5}e^{ay}\left[-\frac{1+\sqrt{5}}{2}
\cos(\sqrt{-b}y)+\sqrt{\frac{5}{-b}}\sin(\sqrt{-b}y)\right]\nonumber}\\
&&+\frac{1}{5}e^{-(1+a)y}\left[\frac{-1+\sqrt{5}}{2}
\cos(\sqrt{5+b}y)-\frac{5+\sqrt{5}}{2\sqrt{-b}}\sin(\sqrt{5+b}y)\right],
\label{g5-22}
\end{eqnarray}
\begin{eqnarray}
\lefteqn{g_{53}(2y)=\frac{1}{5}e^{2y}
+\frac{1}{5}e^{ay}\left[-\frac{1+\sqrt{5}}{2}
\cos(\sqrt{-b}y)-\sqrt{\frac{5}{-b}}\sin(\sqrt{-b}y)\right]\nonumber}\\
&&+\frac{1}{5}e^{-(1+a)y}\left[\frac{-1+\sqrt{5}}{2}
\cos(\sqrt{5+b}y)+\frac{5+\sqrt{5}}{2\sqrt{-b}}\sin(\sqrt{5+b}y)\right],
\label{g5-23}
\end{eqnarray}
\begin{eqnarray}
\lefteqn{g_{54}(2y)=\frac{1}{5}e^{2y}
+\frac{1}{5}e^{ay}\left[\frac{-1+\sqrt{5}}{2}
\cos(\sqrt{-b}y)-\frac{5+\sqrt{5}}{2\sqrt{-b}}
\sin(\sqrt{-b}y)\right]\nonumber}\\
&&+\frac{1}{5}e^{-(1+a)y}\left[-\frac{1+\sqrt{5}}{2}
\cos(\sqrt{5+b}y)-\sqrt{\frac{5}{-b}}\sin(\sqrt{5+b}y)\right].
\label{g5-24}
\end{eqnarray}
The polar 5-dimensional cosexponential functions can be written as
\begin{equation}
g_{5k}(y)=\frac{1}{5}\sum_{l=0}^{4}\exp\left[y\cos\left(\frac{2\pi l}{5}\right)
\right]
\cos\left[y\sin\left(\frac{2\pi l}{5}\right)-\frac{2\pi kl}{5}\right],
k=0,...,4. 
\label{30}
\end{equation}
The graphs of the polar 5-dimensional cosexponential functions are shown in
Fig. 2.

It can be checked that
\begin{equation}
\sum_{k=0}^{4}g_k^2(y)=\frac{1}{5}e^{2y}+\frac{2}{5}e^{(\sqrt{5}-1)y/2}
+\frac{2}{5}e^{-(\sqrt{5}+1)y/2}.
\label{34a}
\end{equation}

The addition theorems for the polar 5-dimensional cosexponential functions are
\begin{eqnarray}
\lefteqn{\begin{array}{c}
g_{50}(y+z)=g_{50}(y)g_{50}(z)
+g_{51}(y)g_{54}(z)+g_{52}(y)g_{53}(z)+g_{53}(y)g_{52}(z)+g_{54}(y)g_{51}(z) 
,\\
g_{51}(y+z)=g_{50}(y)g_{51}(z)+g_{51}(y)g_{50}(z)
+g_{52}(y)g_{54}(z)+g_{53}(y)g_{53}(z)+g_{54}(y)g_{52}(z) ,\\
g_{52}(y+z)=g_{50}(y)g_{52}(z)+g_{51}(y)g_{51}(z)+g_{52}(y)g_{50}(z)
+g_{53}(y)g_{54}(z)+g_{54}(y)g_{53}(z) ,\\
g_{53}(y+z)=g_{50}(y)g_{53}(z)+g_{51}(y)g_{52}(z)
+g_{52}(y)g_{51}(z)+g_{53}(y)g_{50}(z)
+g_{54}(y)g_{54}(z) ,\\
g_{54}(y+z)=g_{50}(y)g_{54}(z)
+g_{51}(y)g_{53}(z)+g_{52}(y)g_{52}(z)
+g_{53}(y)g_{51}(z)+g_{54}(y)g_{50}(z) .\\
\end{array}\nonumber}\\
&&
\label{35a}
\end{eqnarray}
It can be shown that
\begin{eqnarray}
\begin{array}{l}
\{g_{50}(y)+h_1g_{51}(y)+h_2g_{52}(y)+h_3g_{53}(y)+h_4g_{54}(y)\}^l\\
\hspace*{0.5cm}=g_{50}(ly)+h_1g_{51}(ly)+h_2g_{52}(ly)
+h_3g_{53}(ly)+h_4g_{54}(ly),\\
\{g_{50}(y)+h_1g_{53}(y)+h_2g_{51}(y)+h_3g_{54}(y)+h_4g_{52}(y)\}^l\\
\hspace*{0.5cm}=g_{50}(ly)+h_1g_{53}(ly)
+h_2g_{51}(ly)+h_3g_{54}(ly)+h_4g_{52}(ly),\\
\{g_{50}(y)+h_1g_{52}(y)+h_2g_{54}(y)+h_3g_{51}(y)+h_4g_{53}(y)\}^l\\
\hspace*{0.5cm}=g_{50}(ly)+h_1g_{52}(ly)+h_2g_{54}(ly)
+h_3g_{51}(ly)+h_4g_{53}(ly),\\
\{g_{50}(y)+h_1g_{54}(y)+h_2g_{53}(y)+h_3g_{52}(y)+h_4g_{51}(y)\}^l\\
\hspace*{0.5cm}=g_{50}(ly)+h_1g_{54}(ly)+h_2g_{53}(ly)
+h_3g_{52}(ly)+h_4g_{51}(ly).
\end{array}
\label{37b}
\end{eqnarray}

The derivatives of the polar cosexponential functions
are related by
\begin{equation}
\frac{dg_{50}}{du}=g_{54}, \:
\frac{dg_{51}}{du}=g_{50}, \:
\frac{dg_{52}}{du}=g_{51}, \:
\frac{dg_{53}}{du}=g_{52} ,\:
\frac{dg_{54}}{du}=g_{53}. \:
\label{45}
\end{equation}

\section{Exponential and trigonometric forms of polar 5-complex numbers}

The exponential and trigonometric forms of 5-complex
numbers can be expressed with the aid of the hypercomplex bases 
\begin{equation}
\left(
\begin{array}{c}
e_+\\
e_1\\
\tilde e_1\\
e_2\\
\tilde e_2\\
\end{array}\right)
=
\frac{2}{5}\left(
\begin{array}{ccccc}
\frac{1}{2}&\frac{1}{2}&\frac{1}{2}&\frac{1}{2}&\frac{1}{2}\\
1&p&2p^2-1&2p^2-1&p\\
0&q&2pq&-2pq&-q\\
1&2p^2-1&p&p&2p^2-1\\
0&2pq&-q&q&-2pq\\
\end{array}
\right)
\left(
\begin{array}{c}
1\\
h_1\\
h_2\\
h_3\\ 
h_4
\end{array}
\right).
\label{e12}
\end{equation}

The multiplication relations for these bases are
\begin{eqnarray}
\lefteqn{e_+^2=e_+,\;  e_+e_k=0,\; e_+\tilde e_k=0,\; \nonumber}\\ 
&&e_k^2=e_k,\; \tilde e_k^2=-e_k,\; e_k \tilde e_k=\tilde e_k ,\; e_ke_l=0,\;
e_k\tilde e_l=0,\; \tilde e_k\tilde e_l=0,\; k,l=1,2, \;k\not=l.
\label{e12b}
\end{eqnarray}
The bases have the property that
\begin{equation}
e_+ + e_1 +e_2=1.
\label{48b}
\end{equation}
The moduli of the new bases are
\begin{equation}
|e_+|=\frac{1}{\sqrt{5}},\; 
|e_k|=\sqrt{\frac{2}{5}},\; |\tilde e_k|=\sqrt{\frac{2}{5}}, 
\label{e12c}
\end{equation}
for $k=1,2$.

It can be checked that
\begin{eqnarray}
x_0+h_1x_1+h_2x_2+h_3x_3+h_4x_4=
e_+v_+ 
+e_1 v_1+\tilde e_1 \tilde v_1+e_2 v_2+\tilde e_2 \tilde v_2.
\label{e13b}
\end{eqnarray}
The ensemble $e_+, e_1, \tilde e_1, e_2, \tilde e_2$ will be called the
canonical 5-complex base, and Eq. (\ref{e13b}) gives the canonical form of the
5-complex number.

The exponential form of the 5-complex number $u$ is
\begin{eqnarray}
\lefteqn{u=\rho\exp\left\{\frac{1}{5}(h_1+h_2+h_3+h_4)
\ln\frac{\sqrt{2}}{\tan\theta_+}\right.\nonumber}\\
&&\left.+\left[\frac{\sqrt{5}+1}{10}(h_1+h_4)
-\frac{\sqrt{5}-1}{10}(h_2+h_3)\right]
\ln\tan\psi_1+\tilde e_1\phi_1+\tilde e_2\phi_2
\right\},
\label{50b}
\end{eqnarray}
for $0<\theta_+<\pi/2$.

The trigonometric form of the 5-complex number $u$ is
\begin{eqnarray}
\lefteqn{u=d\left(\frac{5}{2}\right)^{1/2}
\left(\frac{1}{\tan^2\theta_+}+1
+\frac{1}{\tan^2\psi_1}\right)^{-1/2}
\left(\frac{e_+\sqrt{2}}{\tan\theta_+}
+e_1+\frac{e_2}{\tan\psi_1}\right)
\exp\left(\tilde e_1\phi_1+\tilde e_2\phi_2\right).\nonumber}\\
&&
\label{52b}
\end{eqnarray}

The modulus $d$ and the amplitude $\rho$ are related by
\begin{eqnarray}
d=\rho \frac{2^{2/5}}{\sqrt{5}}
\left(\tan\theta_+
\tan^2\psi_1\right)^{1/5}
\left(\frac{1}{\tan^2\theta_+}+1
+\frac{1}{\tan^2\psi_1}\right)^{1/2}.
\label{53b}
\end{eqnarray}

\section{Elementary functions of a polar 5-complex variable}

The logarithm and power function exist for $v_+>0$, which means that
$0<\theta_+<\pi/2$, and are given by 
\begin{eqnarray}
\lefteqn{\ln u=\ln\rho+
\frac{1}{5}(h_1+h_2+h_3+h_4)\ln\frac{\sqrt{2}}{\tan\theta_+}\nonumber}\\
&&+\left[\frac{\sqrt{5}+1}{10}(h_1+h_4)-\frac{\sqrt{5}-1}{10}(h_2+h_3)\right]
\ln\tan\psi_1+\tilde e_1\phi_1+\tilde e_2\phi_2,
\label{56b}
\end{eqnarray}
\begin{equation}
u^m=e_+ v_+^m +
\rho_1^m(e_1\cos m\phi_1+\tilde e_1\sin m\phi_1)
+\rho_2^m(e_2\cos m\phi_2+\tilde e_2\sin m\phi_2).
\label{59b}
\end{equation}

The exponential of the 5-complex variable $u$ is
\begin{eqnarray}
e^u=e_+e^{v_+}  
+e^{v_1}\left(e_1 \cos \tilde v_1+\tilde e_1 \sin\tilde v_1\right)
+e^{v_2}\left(e_2 \cos \tilde v_2+\tilde e_2 \sin\tilde v_2\right).
\label{73b}
\end{eqnarray}
The trigonometric functions of the
5-complex variable $u$ are
\begin{equation}
\cos u=e_+\cos v_+  
+\sum_{k=1}^{2}\left(e_k \cos v_k\cosh \tilde v_k
-\tilde e_k \sin v_k\sinh\tilde v_k\right),
\label{74c}
\end{equation}
\begin{equation}
\sin u=e_+\sin v_+  
+\sum_{k=1}^{2}\left(e_k \sin v_k\cosh \tilde v_k
+\tilde e_k \cos v_k\sinh\tilde v_k\right).
\label{74d}
\end{equation}

The hyperbolic functions of the
5-complex variable $u$ are
\begin{equation}
\cosh u=e_+\cosh v_+  
+\sum_{k=1}^{2}\left(e_k \cosh v_k\cos \tilde v_k
+\tilde e_k \sinh v_k\sin\tilde v_k\right),
\label{75c}
\end{equation}
\begin{equation}
\sinh u=e_+\sinh v_+  
+\sum_{k=1}^{2}\left(e_k \sinh v_k\cos \tilde v_k
+\tilde e_k \cosh v_k\sin\tilde v_k\right).
\label{75d}
\end{equation}

\section{Power series of 5-complex numbers}

A power series of the 5-complex variable $u$ is a series of the form
\begin{equation}
a_0+a_1 u + a_2 u^2+\cdots +a_l u^l+\cdots .
\label{83}
\end{equation}
Since
\begin{equation}
|au^l|\leq 5^{l/2} |a| |u|^l ,
\label{82}
\end{equation}
the series is absolutely convergent for 
\begin{equation}
|u|<c,
\label{86}
\end{equation}
where 
\begin{equation}
c=\lim_{l\rightarrow\infty} \frac{|a_l|}{\sqrt{5}|a_{l+1}|} .
\label{87}
\end{equation}

If $a_l=\sum_{p=0}^{4}h_p a_{lp}$, where $h_0=1$, and
\begin{equation}
A_{l+}=\sum_{p=0}^{4}a_{lp},
\label{88a}
\end{equation}
\begin{equation}
A_{lk}=\sum_{p=0}^{4}a_{lp}\cos\left(\frac{2\pi kp}{5}\right),
\label{88b}
\end{equation}
\begin{equation}
\tilde A_{lk}=\sum_{p=0}^{4}a_{lp}\sin\left(\frac{2\pi kp}{5}\right),
\label{88c}
\end{equation}
for $k=1,2$,
the series (\ref{83}) can be written as
\begin{equation}
\sum_{l=0}^\infty \left[
e_+A_{l+}v_+^l+\sum_{k=1}^{2}
(e_k A_{lk}+\tilde e_k\tilde A_{lk})(e_k v_k+\tilde e_k\tilde v_k)^l 
\right].
\label{89b}
\end{equation}
The series in Eq. (\ref{83}) is absolutely convergent for 
\begin{equation}
|v_+|<c_+,\:
\rho_k<c_k, k=1,2,
\label{90}
\end{equation}
where 
\begin{equation}
c_+=\lim_{l\rightarrow\infty} \frac{|A_{l+}|}{|A_{l+1,+}|} ,\:
c_k=\lim_{l\rightarrow\infty} \frac
{\left(A_{lk}^2+\tilde A_{lk}^2\right)^{1/2}}
{\left(A_{l+1,k}^2+\tilde A_{l+1,k}^2\right)^{1/2}} .
\label{91}
\end{equation}

\section{Analytic functions of a polar 5-compex variable}

The expansion of an analytic function $f(u)$ around $u=u_0$ is
\begin{equation}
f(u)=\sum_{k=0}^\infty \frac{1}{k!} f^{(k)}(u_0)(u-u_0)^k .
\label{h91d}
\end{equation}
Since the limit $f^\prime (u_0)=\lim_{u\rightarrow u_0}\{f(u)-f(u_0)\}/(u-u_0)$

is independent of the
direction in space along which $u$ is approaching $u_0$, the function $f(u)$ 
is said to be analytic, analogously to the case of functions of regular complex
variables. \cite{3} 
If $f(u)=\sum_{k=0}^{4}h_kP_k(x_0,x_1,x_2,x_3,x_{4})$, then
\begin{equation}
\frac{\partial P_0}{\partial x_0} 
=\frac{\partial P_1}{\partial x_1} 
=\frac{\partial P_2}{\partial x_2} 
=\frac{\partial P_3}{\partial x_3}
=\frac{\partial P_4}{\partial x_4}, 
\label{h95a}
\end{equation}
\begin{equation}
\frac{\partial P_1}{\partial x_0} 
=\frac{\partial P_2}{\partial x_1} 
=\frac{\partial P_3}{\partial x_2} 
=\frac{\partial P_4}{\partial x_3}
=\frac{\partial P_0}{\partial x_4}, 
\label{h95b}
\end{equation}
\begin{equation}
\frac{\partial P_2}{\partial x_0} 
=\frac{\partial P_3}{\partial x_1} 
=\frac{\partial P_4}{\partial x_2} 
=\frac{\partial P_0}{\partial x_3}
=\frac{\partial P_1}{\partial x_4}, 
\label{h95c}
\end{equation}
\begin{equation}
\frac{\partial P_3}{\partial x_0} 
=\frac{\partial P_4}{\partial x_1} 
=\frac{\partial P_0}{\partial x_2} 
=\frac{\partial P_1}{\partial x_3}
=\frac{\partial P_2}{\partial x_4}, 
\label{h95d}
\end{equation}
\begin{equation}
\frac{\partial P_4}{\partial x_0} 
=\frac{\partial P_0}{\partial x_1} 
=\frac{\partial P_1}{\partial x_2} 
=\frac{\partial P_2}{\partial x_3}
=\frac{\partial P_3}{\partial x_4}, 
\label{h95e}
\end{equation}
and
\begin{eqnarray}
\lefteqn{\frac{\partial^2 P_k}{\partial x_0\partial x_l}
=\frac{\partial^2 P_k}{\partial x_1\partial x_{l-1}}
=\cdots=
\frac{\partial^2 P_k}{\partial x_{[l/2]}\partial x_{l-[l/2]}}}\nonumber\\
&&=\frac{\partial^2 P_k}{\partial x_{l+1}\partial x_{4}}
=\frac{\partial^2 P_k}{\partial x_{l+2}\partial x_{3}}
=\cdots
=\frac{\partial^2 P_k}{\partial x_{l+1+[(3-l)/2]}
\partial x_{4-[(3-l)/2]}} ,
\label{96}
\end{eqnarray}
for $k,l=0,...,4$.
In Eq. (\ref{96}), $[a]$ denotes the integer part of $a$,
defined as $[a]\leq a<[a]+1$. 
In this work, brackets larger than the regular brackets
$[\;]$ do not have the meaning of integer part.

\section{Integrals of polar 5-complex functions}

If $f(u)$ is an analytic 5-complex function,
then
\begin{equation}
\oint_\Gamma \frac{f(u)du}{u-u_0}=
2\pi f(u_0)\left\{\tilde e_1 
\;{\rm int}(u_{0\xi_1\eta_1},\Gamma_{\xi_1\eta_1})+
\tilde e_2 
\;{\rm int}(u_{0\xi_2\eta_2},\Gamma_{\xi_2\eta_2})\right\} ,
\label{120}
\end{equation}
where
\begin{equation}
{\rm int}(M,C)=\left\{
\begin{array}{l}
1 \;\:{\rm if} \;\:M \;\:{\rm is \;\:an \;\:interior \;\:point \;\:of} 
\;\:C ,\\ 
0 \;\:{\rm if} \;\:M \;\:{\rm is \;\:exterior \;\:to}\:\; C ,\\
\end{array}\right.,
\label{118}
\end{equation}
and $u_{0\xi_k\eta_k}$, $\Gamma_{\xi_k\eta_k}$ are respectively the
projections of the pole $u_0$ and of 
the loop $\Gamma$ on the plane defined by the axes $\xi_k$ and $\eta_k$,
$k=1,2$.

\section{Factorization of polar 5-complex polynomials}

A polynomial of degree $m$ of the 5-complex variable $u$ has the form
\begin{equation}
P_m(u)=u^m+a_1 u^{m-1}+\cdots+a_{m-1} u +a_m ,
\label{125}
\end{equation}
where $a_l$, for $l=1,...,m$, are 5-complex constants.
If $a_l=\sum_{p=0}^{4}h_p a_{lp}$, and with the
notations of Eqs. (\ref{88a})-(\ref{88c}) applied for $l= 1, \cdots, m$, the
polynomial $P_m(u)$ can be written as 
\begin{eqnarray}
\lefteqn{P_m= 
e_+\left(v_+^m +\sum_{l=1}^{m}A_{l+}v_+^{m-l} \right) \nonumber}\\
&&+\sum_{k=1}^{2}
\left[(e_k v_k+\tilde e_k\tilde v_k)^m+
\sum_{l=1}^m(e_k A_{lk}+\tilde e_k\tilde A_{lk})
(e_k v_k+\tilde e_k\tilde v_k)^{m-l} 
\right].
\label{126b}
\end{eqnarray}

The polynomial $P_m(u)$ can be written, as 
\begin{eqnarray}
P_m(u)=\prod_{p=1}^m (u-u_p) ,
\label{128c}
\end{eqnarray}
where
\begin{eqnarray}
u_p=e_+ v_{p+}
+\left(e_1 v_{1p}+\tilde e_1\tilde v_{1p}\right)
+\left(e_2 v_{2p}+\tilde e_2\tilde v_{2p}\right), p=1,...,m.
\label{128e}
\end{eqnarray}
The quantities $v_{p+}$,   
$e_k v_{kp}+\tilde e_k\tilde v_{kp}$,
$p=1,...,m, k=1,2$,
are the roots of the corresponding polynomial in Eq. (\ref{126b}). The roots
$v_{p+}$ appear in complex-conjugate pairs, and  
$v_{kp}, \tilde v_{kp}$ are real numbers.
Since all these roots may be ordered arbitrarily, the polynomial $P_m(u)$ can
be written in many different ways as a product of linear factors.

If $P(u)=u^2-1$, the degree is $m=2$, the coefficients of the polynomial are
$a_1=0, a_2=-1$, the coefficients defined in Eqs. (\ref{88a})-(\ref{88c})
are $A_{2+}=-1, A_{21}=-1, \tilde A_{21}=0,
A_{22}=-1, \tilde A_{22}=0$. The expression of $P(u)$, Eq. (\ref{126b}), is  
$v_+^2-e_++(e_1v_1+\tilde e_1\tilde v_1)^2-e_1+
(e_2v_2+\tilde e_2\tilde v_2)^2-e_2 $. 
The factorization of $P(u)$, Eq. (\ref{128c}), is
$P(u)=(u-u_1)(u-u_2)$, where the roots are
$u_1=\pm e_+\pm e_1\pm  e_2, u_2=-u_1$. If $e_+, e_1, e_2$ 
are expressed with the aid of Eq. (\ref{e12}) in terms of $h_1, h_2, h_3,
h_4$, the factorizations of $P(u)$ are obtained as
\begin{eqnarray}
\lefteqn{\begin{array}{l}
u^2-1=(u+1)(u-1),\\
u^2-1=\left[u+\frac{1}{5}+\frac{\sqrt{5}+1}{5}(h_1+h_4)
-\frac{\sqrt{5}-1}{5}(h_2+h_3)\right]
\left[u-\frac{1}{5}-\frac{\sqrt{5}+1}{5}(h_1+h_4)
+\frac{\sqrt{5}-1}{5}(h_2+h_3)\right],\\
u^2-1=\left[u+\frac{1}{5}-\frac{\sqrt{5}-1}{5}(h_1+h_4)
+\frac{\sqrt{5}+1}{5}(h_2+h_3)\right]
\left[u-\frac{1}{5}+\frac{\sqrt{5}-1}{5}(h_1+h_4)
-\frac{\sqrt{5}+1}{5}(h_2+h_3)\right],\\
u^2-1=\left[u+\frac{3}{5}-\frac{2}{5}(h_1+h_2+h_3+h_4)\right]
\left[u-\frac{3}{5}+\frac{2}{5}(h_1+h_2+h_3+h_4)\right].
\end{array}\nonumber}\\
&&
\end{eqnarray}
It can be checked that 
$(\pm e_+\pm e_1\pm e_2)^2=e_++e_1+e_2=1$.

\section{Representation of polar 5-complex numbers by irreducible matrices}

If the unitary matrix which can be obtained from the expression, 
Eq. (\ref{12}), of the variables $\xi_+, \xi_1, \eta_1, \xi_k, \eta_k$ in terms
of $x_0, x_1, x_2, x_3, x_4$ is called $T$,
the irreducible representation \cite{4} of the hypercomplex number $u$ is
\begin{equation}
T U T^{-1}=\left(
\begin{array}{ccc}
v_+     &     0     &     0      \\
0       &     V_1   &     0      \\
0       &     0     &     V_2    \\
\end{array}
\right),
\label{129b}
\end{equation}
where $U$ is the matrix in Eq. (\ref{24b}),
and $V_k$ are the matrices
\begin{equation}
V_k=\left(
\begin{array}{cc}
v_k           &     \tilde v_k   \\
-\tilde v_k   &     v_k          \\
\end{array}\right),\;\; k=1,2.
\label{130}
\end{equation}

\section{Conclusions}

The operations of addition and multiplication of the 5-complex numbers
introduced in this 
work have a geometric interpretation based on the amplitude $\rho$,
the modulus $d$ and the polar, planar and azimuthal angles $\theta_+, 
\psi_1, \phi_1,\phi_2$. 
If $x_0+x_1+x_2+x_3+x_4>0$
the 5-complex numbers can be written in exponential and
trigonometric forms with the aid of the modulus, amplitude and the angular
variables. 
The 5-complex functions defined by series of powers are analytic, and 
the partial derivatives of the components of the 5-complex functions are
closely related. The integrals of 5-complex functions are independent of path
in regions where the functions are regular. The fact that the exponential form
of the 5-complex numbers depends on the cyclic variables $\phi_k$
leads to the 
concept of pole and residue for integrals on closed paths. The polynomials of
5-complex variables can be written as products of linear or quadratic
factors.

\newpage

FIGURE CAPTIONS\\

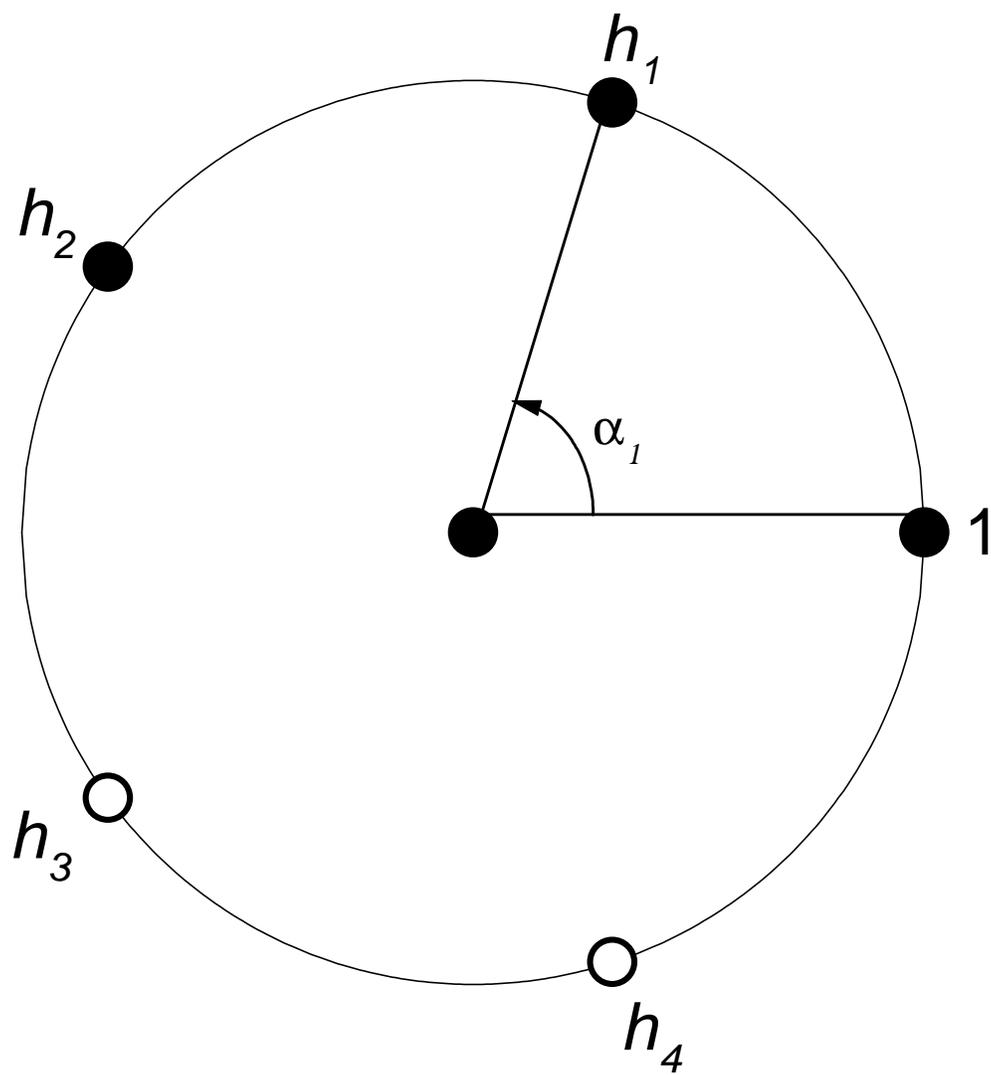
Fig. 1. Representation of the polar 
hypercomplex bases $1,h_1,h_2,h_3,h_4$
by points on a circle at the angles $\alpha_k=2\pi k/5$.
The product $h_j h_k$ will be represented by the point of the circle at the
angle $2\pi (j+k)/5$, $i,k=0,1,...,4$, where $h_0=1$. If $2\pi\leq 
2\pi (j+k)/5\leq 4\pi$, the point represents the basis
$h_l$ of angle $\alpha_l=2\pi(j+k)/5-2\pi$.\\

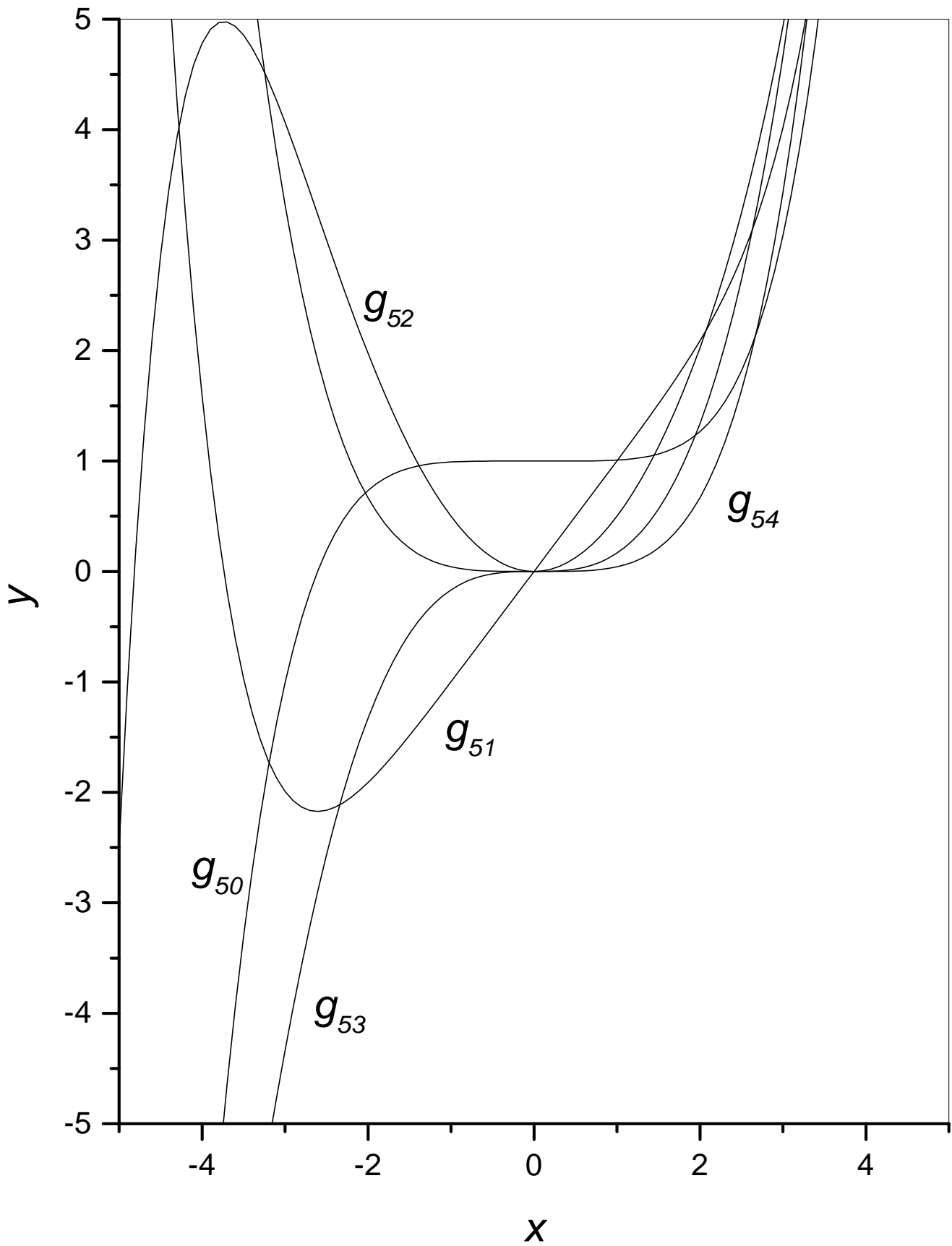
Fig. 2. Polar cosexponential functions 
$g_{50}, g_{51},g_{52}, g_{53},g_{54}$.

\end{document}